\documentclass[12pt]{amsart}
\usepackage{amssymb}
\usepackage{amscd}
\usepackage{pstricks}

\newcommand{\Q}{\mathbb{Q}}
\newcommand{\Z}{\mathbb{Z}}
\newcommand{\N}{\mathbb{N}}
\newcommand{\R}{\mathbb{R}}
\newcommand{\C}{\mathbb{C}}
\newcommand{\Ss}{\mathbb{S}}
\newcommand{\Aa}{\mathbb{A}}
\newcommand{\Pp}{\mathbb{P}}
\newcommand{\bbone}{\mathbf{1}}

\newcommand{\CalA}{\mathcal{A}}

\newcommand{\CalP}{\mathcal{P}}
\newcommand{\CalG}{\mathcal{G}}
\newcommand{\CalT}{\mathcal{T}}
\newcommand{\CalM}{\mathcal{M}}

\newcommand{\tH}{\widetilde{H}}
\newcommand{\WH}{W\! H}
\newcommand{\YBn}{\overline{Y}_{B_n}}
\newcommand{\YDn}{\overline{Y}_{D_n}}

\newcommand{\Exp}{\mathrm{Exp}}
\newcommand{\Sinh}{\mathrm{Sinh}}
\newcommand{\Cosh}{\mathrm{Cosh}}
\newcommand{\Sech}{\mathrm{Sech}}
\newcommand{\sech}{\mathrm{sech}}

\newcommand{\Arcsinh}{\mathrm{Arcsinh}}
\newcommand{\arcsinh}{\mathrm{arcsinh}}
\newcommand{\Ind}{\mathrm{Ind}}

\newcommand{\Hom}{\mathrm{Hom}}

\newcommand{\tr}{\mathrm{tr}}

\newcommand{\rk}{\mathrm{rk}}
\newcommand{\ch}{\mathrm{ch}}

\newcommand{\rmprod}{\mathrm{prod}}

\newcommand{\orient}{\mathrm{or}}

\newcommand{\lng}{\langle}
\newcommand{\rng}{\rangle}
\newcommand{\sign}{\sim}

\numberwithin{equation}{section}
\newtheorem{theorem}{Theorem}[section]
\newtheorem{corollary}[theorem]{Corollary}
\newtheorem{lemma}[theorem]{Lemma}
\newtheorem{proposition}[theorem]{Proposition}

\theoremstyle{definition}
\newtheorem{definition}[theorem]{Definition}

\newtheorem{example}[theorem]{Example}

\newcommand{\bdf}{\begin{definition}}
\newcommand{\edf}{\end{definition}\noindent}
\newcommand{\bex}{\begin{example}}
\newcommand{\eex}{\end{example}\noindent}
\newcommand{\bpr}{\begin{proposition}}
\newcommand{\epr}{\end{proposition}}
\newcommand{\blm}{\begin{lemma}}
\newcommand{\elm}{\end{lemma}}
\newcommand{\bth}{\begin{theorem}}
\renewcommand{\eth}{\end{theorem}}
\newcommand{\bpf}{\begin{proof}}
\newcommand{\epf}{\end{proof}\noindent}
\newcommand{\bcr}{\begin{corollary}}
\newcommand{\ecr}{\end{corollary}\noindent}
\newcommand{\beq}{\begin{equation}}
\newcommand{\eeq}{\end{equation}}
\newcommand{\bes}{\begin{equation*}}
\newcommand{\ees}{\end{equation*}}
\newcommand{\ben}{\begin{enumerate}}
\newcommand{\een}{\end{enumerate}}
\begin{document}
\title[Real De~Concini--Procesi models of Coxeter type]
{The cohomology of real De~Concini--Procesi models of Coxeter type}
\author{Anthony Henderson}
\address{School of Mathematics and Statistics,
University of Sydney, NSW 2006, AUSTRALIA}
\email{anthonyh@maths.usyd.edu.au}
\thanks{This work was supported by Australian Research Council grant DP0344185}
\author{Eric Rains}
\address{Department of Mathematics, California Institute of Technology,
Pasadena, CA 91125, USA}
\email{rains@caltech.edu}
\begin{abstract}
We study the rational cohomology groups of the real
De~Concini--Procesi model corresponding to a finite Coxeter group,
generalizing the type-$A$ case of the moduli space of stable
genus $0$ curves with marked points. We compute the Betti numbers
in the exceptional types, and give formulae for them in types
$B$ and $D$. We give a generating-function formula for the characters of
the representations of a Coxeter group of type $B$ on the
rational cohomology groups of the corresponding
real De~Concini--Procesi model, and
deduce the multiplicities of one-dimensional
characters in the representations, and a formula for the Euler
character. We also give a moduli space interpretation of this
type-$B$ variety, and hence show that the action of the Coxeter group
extends to a slightly larger group.
\end{abstract}
\maketitle
\section*{Introduction}
In \cite{etingofetal}, the second author and 
his collaborators studied the rational cohomology ring of 
$\overline{\CalM_{0,n}}(\R)$, the manifold of real points of the moduli
space of stable genus $0$ curves with $n$ marked points. 
In \cite{rainstypea},
the second author described this ring as a representation of $S_n$,
giving a formula for the graded character. The crux of both papers
was a close connection between the cohomology of
$\overline{\CalM_{0,n}}(\R)$ and the homology of the poset
of partitions of $\{1,2,\cdots,n-1\}$ into parts of odd size.
One explanation for this connection is
that $\overline{\CalM_{0,n}}(\R)$ can be regarded as the real
De~Concini--Procesi model of the hyperplane arrangement of type
$A_{n-2}$; in this context the relation with poset homology
was generalized to an arbitrary subspace arrangement in \cite{rainshomology}.
It was then natural to try to extend the results
of \cite{etingofetal} and \cite{rainstypea} to the
other finite Coxeter types, and that is the goal of the present work.

For any irreducible finite Coxeter group $W$ with reflection
representation $V$ over $\R$, the real 
De~Concini--Procesi model $\overline{Y}_W(\R)$ is defined 
(see Section 1) to be the closure of the image of the map
\[ M_W \to\prod_{W'\in\Pi_W^{\mathrm{irr}}}\Pp(V/V^{W'}), \]
where $M_W$ is the complement of the reflecting hyperplanes and
$\Pi_W^{\mathrm{irr}}$ denotes the set of irreducible parabolic subgroups
of $W$. This is a nonsingular projective real variety. If
$W=S_{n-1}$, $\overline{Y}_W(\R)$ is isomorphic to 
$\overline{\CalM_{0,n}}(\R)$, which has been studied in \cite{devadoss}
and \cite{davisetal} as well as the papers mentioned above.
The general Coxeter case has been considered in
\cite{davisetal}, \cite{gaiffi}, and \cite{armstrongetal},
from varying points of view and with different emphases from the present
paper, which concentrates on the rational cohomology groups
$H^i(\overline{Y}_W(\R);\Q)$. The cohomology groups of the
corresponding complex varieties have also been studied --
see, for instance, \cite{ginzburgkapranov} for
$\overline{\CalM_{0,n}}(\C)$, \cite{yuzvinsky} for
the classical types, and \cite{hendersonwreath} for type $B$ --
but these behave very differently from the real case.

In Section 1 of this paper we recall the main result of \cite{rainshomology},
expressing the cohomology of a real De~Concini--Procesi model in terms
of the homology of a corresponding poset: in the case of 
$\overline{Y}_W(\R)$, the poset $\Pi_W^{(2)}$ consists of those
parabolic subgroups of $W$ whose irreducible components all have even rank.
We prove in Theorem \ref{posetthm}
that this poset is Cohen--Macaulay,
and deduce in Theorem \ref{rainscoxeterthm}
that $H^i(\overline{Y}_W(\R);\Q)$ coincides with the $i$th
Whitney homology $\WH_i(\Pi_W^{(2)})$, up to some twisting
by sign representations. We also prove a result (Corollary \ref{typedcor})
reducing most questions in type $D$ to type $B$.

In Section 2 we compute the Poincar\'e polynomial 
and Euler characteristic of $\overline{Y}_W(\R)$
(these equal the characteristic polynomial and Euler characteristic
of the poset $\Pi_W^{(2)}$). This section is brief, since 
type $A$ was done in \cite{etingofetal}, type $B$ was effectively
done in \cite{hendersondowling} (though the Poincar\'e polynomial
formula is not as explicit as in type $A$), type $D$ reduces to
type $B$, and the other types can be calculated directly.
(The Euler characteristics in the classical types were
calculated in \cite{devadoss}, \cite{gaiffi}, and \cite{armstrongetal}).

In Theorem \ref{mainthm} we prove a plethystic generating-function
formula
for the characters of the representations of $W(B_n)$ on
$H^i(\YBn(\R);\Q)$. This is analogous to the type-$A$ result
in \cite{rainstypea}, but does not give as explicit a formula
for the graded character of a particular group element.
However, we are able to deduce
the multiplicities of the one-dimensional characters
(Proposition \ref{onedimcharprop}), as well as a formula for the Euler
character (Proposition \ref{eulercharprop}). 

In Theorem \ref{modulithm} we give a moduli space intepretation of the
variety $\overline{Y}_{B_n}$: it can
be embedded as a closed subvariety in $\overline{\CalM_{0,2n+2}}$,
namely the one parametrizing stable curves possessing an involution
which fixes the last two marked points, interchanges the $j$th and $(n+j)$th
for $1\leq j\leq n$, and does not fix any component pointwise.
Swapping the labels of the last two marked points induces an involution
$\sigma$ of $\overline{Y}_{B_n}$ which commutes with the $W(B_n)$-action,
and we give some partial results about the action of this
involution on the cohomology of $\overline{Y}_{B_n}(\R)$.
\section{Real De~Concini--Procesi models of Coxeter type}
Let $V$ be a finite-dimensional real vector space, and
$\CalG$ a building set in the dual space $V^*$. (See \cite{wonderful}
for the definition of building set, and of the building set(s)
associated to a subspace arrangement in $V$.) 
For any $G\in\CalG$, let $G^\perp$ denote
the orthogonal subspace of $V$, and let $M_\CalG$ be the
complement $V\setminus\bigcup_{G\in\CalG}G^\perp$.
The associated (compact) real De~Concini--Procesi model 
$\overline{Y}_\CalG(\R)$
is the closure of the image of the map
\[ M_\CalG \to\prod_{G\in\CalG}\Pp(V/G^\perp) \]
induced by the natural maps $M_\CalG\to(V/G^\perp)\setminus\{0\}\to
\Pp(V/G^\perp)$.
Here it makes no difference whether the closure is taken in the Zariski
or in the usual topology, because
Gaiffi's result \cite[Theorem 4.1]{gaiffi} shows that the topological closure
coincides with the set of real points of the
variety $\overline{Y}_\CalG$ defined using Zariski closure.
It was shown in \cite{wonderful} that $\overline{Y}_\CalG$
is a nonsingular projective variety, and hence
$\overline{Y}_\CalG(\R)$ is a
compact smooth manifold, which is easily seen to be connected.

Let $\Pi_\CalG$ denote the lattice of all subspaces of $V^*$
which can be written as direct sums of elements of $\CalG$, and
write $\Pi_\CalG^{(2)}$ for the sub-poset consisting of
direct sums of even-dimensional elements of $\CalG$.
In \cite[Theorem 3.7]{rainshomology} the second author expressed
the homology groups of $\overline{Y}_\CalG(\R)$
(modulo $2$-torsion) in terms of the poset cohomology groups of
$\Pi_\CalG^{(2)}$.
The rational cohomology groups
of $\overline{Y}_\CalG$ can be similarly expressed in terms
of the rational homology groups of $\Pi_\CalG^{(2)}$.
For any $A\in\Pi_\CalG^{(2)}\setminus\{0\}$,
let $H_i(0,A)$ denote the reduced poset homology
$\tH_{i-2}((0,A);\Q)$ of the open interval $(0,A)$, and define
$H_i(0,0)$ to be $\Q$ if $i=0$ and $0$ otherwise. Also
let $\orient(A)$ denote $\tH_{\dim A -1}(\Pp A;\Q)$; since $\dim A$
is even, this is a one-dimensional
vector space on which $GL(A)$ acts via the sign of the determinant.
\bth \label{rainsthm}
For all $i$, there is an isomorphism of $\Q$-vector spaces
\[ H^i(\overline{Y}_\CalG(\R);\Q)\cong 
\bigoplus_{A\in\Pi_\CalG^{(2)}} H_{\dim A - i}(0,A)\otimes\orient(A) \]
which is equivariant for the subgroup of $GL(V)$ preserving $\CalG$.
\eth
\bpf
Apply $\Hom_\Z(-,\Q)$ to the isomorphism in 
\cite[Theorem 3.7]{rainshomology}. Equivariance follows from
the naturality shown there.
\epf
Observe the similarity between this result and that of
Goresky and MacPherson (equivariantly, Sundaram and Welker -- see
\cite[Section 5.4]{wachs})
relating the cohomology of $M_\CalG$
to the homology of the lattice $\Pi_\CalG$.

In this paper we consider building sets associated to
hyperplane arrangements of finite Coxeter type. Let $\CalA$ be such
a hyperplane arrangement in $V$, which we may as well assume to be
essential (i.e.\ the intersection of all the hyperplanes is $0$),
and let $W$ be the corresponding finite Coxeter group. 
We write $n$ for $\rk(W)=\dim V$.
For any hyperplane $H\in\CalA$, let $\alpha_H\in V^*$
be a linear form such that $H=\ker\alpha_H$. Let $\CalG$ be the
minimal building set associated to $\CalA$, namely the 
set of all subspaces $G\subseteq V^*$ which are spanned by some subset
of the $\alpha_H$'s, and cannot be written as a nontrivial direct sum
$G=G_1\oplus G_2$ in such a way that every $\alpha_H\in G$ is contained
in either $G_1$ or $G_2$. Then $\Pi_\CalG$ is the lattice
of all subspaces spanned by a subset of the $\alpha_H$'s.

It is well known that $\Pi_\CalG$ is
isomorphic to the lattice $\Pi_W$ of parabolic subgroups of $W$,
via the bijection sending a parabolic subgroup $W'$ to the subspace
$(V^{W'})^\perp\subseteq V^*$, whose dimension is $\rk(W')$. 
(Note that $\Pi_W$ contains all parabolic subgroups, not just the
standard parabolic subgroups relative to some chosen Coxeter system.)
Under this bijection 
$\CalG$ corresponds to the set $\Pi_W^{\mathrm{irr}}$ of 
irreducible parabolic subgroups, and
the sub-poset $\Pi_\CalG^{(2)}$ corresponds to the poset $\Pi_W^{(2)}$ of
parabolic subgroups of $W$ all of whose irreducible components have even 
rank. The minimum element of $\Pi_W^{(2)}$ will be written $\hat{0}$,
although it is just the trivial subgroup; if all irreducible components
of $W$ have even rank, there is also a maximum element in $\Pi_W^{(2)}$, namely $W$.

Thus we can rephrase the definition of 
$\overline{Y}_W=\overline{Y}_\CalG$ as follows: it is the closure
of the image of the map
\beq \label{modeldefeqn} 
M_W \to\prod_{W'\in\Pi_W^{\mathrm{irr}}}\Pp(V/V^{W'}),
\eeq
where $M_W=V\setminus\bigcup_{H\in\CalA}H$ is the hyperplane complement.
The normalizer $N_{GL(V)}(W)$ preserves $\CalG$ and acts on
$\overline{Y}_W$ by variety automorphisms. The manifold
$\overline{Y}_W(\R)$ is closely related to the minimal blow-up
of the Coxeter complex, studied in \cite{davisetal} and \cite{armstrongetal}
(among other papers);
the latter could be obtained by including in the codomain of
\eqref{modeldefeqn} the factor $S(V)$, the unit sphere in $V$.
(If $W$ is irreducible, the minimal blow-up of the Coxeter complex is
a double cover of $\overline{Y}_W(\R)$.)
Our aim here is to describe the rational cohomology
of $\overline{Y}_W(\R)$, using Theorem \ref{rainsthm}.

It is clear that if $W\cong W_1\times W_2$ is reducible, then
$\Pi_W\cong\Pi_{W_1}\times\Pi_{W_2}$ and
$\overline{Y}_W\cong 
\overline{Y}_{W_1}\times\overline{Y}_{W_2}$, so 
we can restrict attention to the case that
$W$ is irreducible. In this case 
the image of $M_W\to\overline{Y}_W$ is isomorphic to the image $\Pp M_W$ of
$M_W\to\Pp(V)$, which implies that 
$\dim \overline{Y}_W=n-1$.
Note that if $n=1$, then $\overline{Y}_W$ is a point, and if
$n=2$, then $\overline{Y}_W$ is a projective line, so the
interest lies in rank $\geq 3$. When speaking
of a particular type $X_n$, we will write $W$ as $W(X_n)$, $\Pi_W$ as $\Pi_{X_n}$,
$\overline{Y}_W$ as $\overline{Y}_{X_n}$, and so forth.

The crucial question is what special properties the poset
$\Pi_W^{(2)}$ has that a general $\Pi_\CalG^{(2)}$ does not.
We first consider the classical types.
\bex
Suppose $W=W(A_n)\cong S_{n+1}$. Then $\Pi_{A_n}$ is isomorphic to
the partition lattice $\Pi_{n+1}$, and $\Pi_{A_n}^{(2)}$ corresponds
to the poset $\Pi_{n+1}^{1\, \mathrm{mod}\, 2}$ of partitions of
$\{1,\cdots,n+1\}$ into parts of odd size. This poset was shown to be
Cohen--Macaulay by Bj\"orner, and its M\"obius function was studied in \cite{chr}.
\eex 
\bex
Suppose $W=W(B_n)\cong \{\pm 1\}\wr S_n$. Then $\Pi_{B_n}$ is isomorphic
to the Dowling lattice $Q_n(\{\pm 1\})$ (also known as
the signed partition lattice).
Following \cite{hendersondowling}, this consists of pairs
$(J,\pi)$ where $J=\{\pm 1\}\times I$ for some subset 
$I\subseteq\{1,\cdots,n\}$ and $\pi$ is a partition of
$(\{\pm 1\}\times\{1,\cdots,n\})\setminus J$ whose parts are interchanged
in pairs by the action of $\{\pm 1\}$. The sub-poset 
$\Pi_{B_n}^{(2)}$ corresponds to the sub-poset of $Q_n(\{\pm 1\})$ defined
by the conditions $\frac{|J|}{2}\equiv 0$ mod $2$ and $|K|\equiv 1$ mod $2$
for all $K\in\pi$; this is the sub-poset denoted 
$Q_n^{1\text{ mod }2}(\{\pm 1\})$ in 
\cite[Definition 1.3]{hendersondowling}, except that when $n$ is odd
we must remove the maximum element artifically added to make the latter 
poset bounded. So \cite[Proposition 1.4]{hendersondowling} shows that
$\Pi_{B_n}^{(2)}$ is Cohen--Macaulay.
\eex
\bex
Suppose $W=W(D_n)$ for $n\geq 4$. There is an obvious
poset embedding $\Pi_{D_n}^{(2)}\hookrightarrow\Pi_{B_n}^{(2)}$
which replaces any $D_\ell$ factor of a parabolic subgroup by
a $B_\ell$ factor. The image of this embedding is the sub-poset
defined by the condition $\frac{|J|}{2}\neq 2$ (because $D_2=A_1\times A_1$
is reducible). One could prove Cohen--Macaulayness of 
this sub-poset by an argument similar
to \cite[Proposition 1.4]{hendersondowling}, although
there is a slight complication: the principal upper order ideals
corresponding to pairs $(\emptyset,\pi)$ are not in general isomorphic to
a smaller $\Pi_{D_\ell}^{(2)}$, but rather to posets associated to
the hyperplane 
arrangements interpolating between types $B$ and $D$ (i.e.\ those denoted
$\CalA_\ell^k$ in \cite[Proposition 6.82]{orlikterao}). 
\eex
We have the following type-independent result (see
\cite{wachs} for the definitions involved). Note that finiteness
of $W$ is necessary only for (3).
\bth \label{posetthm}
Let $W$ be a finite Coxeter group.
\ben
\item $\Pi_W^{(2)}$ is pure with rank function
$\rk_{\Pi_W^{(2)}}(W')=\frac{\rk(W')}{2}$.
\item If $W$ is irreducible of even rank, then $\Pi_W^{(2)}$ is semimodular;
if $W$ is irreducible of odd rank, then $\Pi_W^{(2)}\cup\{\hat{1}\}$
is semimodular.
\item $\Pi_W^{(2)}$ is Cohen--Macaulay.
\een
\eth
\bpf
Let $I$ be the set of vertices of the Coxeter graph of $W$.
For any $J\subseteq I$, let $W_J$ denote the corresponding standard parabolic
subgroup. It is well known that
every chain of parabolic subgroups of $W$ is $W$-conjugate
to a chain of standard parabolic subgroups.

To prove (1), it suffices to show that if $W'\subset W''$ is a cover
relation in $\Pi_W^{(2)}$, then $\rk(W'')=\rk(W')+2$. Conjugating by
a suitable element of $W$ if necessary, we can assume that $W'=W_J$
and $W''=W_K$ are standard parabolic subgroups, where $J\subset K$. 
Regarding $J$ and $K$
as subgraphs of the Coxeter graph of $W$, we know by definition of
$\Pi_W^{(2)}$ that every connected component of $K$ is even (i.e.\ has an even
number of vertices), and similarly for $J$. Moreover, this evenness property
does not hold for any subgraph $J'$ such that $J\subset J'\subset K$.
It is easy to see that this forces $|K\setminus J|=2$ as required.

In part (2), we assume that $W$ is irreducible, so the Coxeter graph
is connected. We must show that if
$W',W''\in \Pi_W^{(2)}$ are different covers of a single
element $W^0$, then either there is another element of $\Pi_W^{(2)}$
which covers $W'$ and $W''$, or both $W'$ and $W''$ are maximal in
$\Pi_W^{(2)}$ (the latter case being possible only if $\rk(W)$ is odd). 
By part (1), the assumption
implies that $\rk(W')=\rk(W'')=\rk(W^0)+2$. Hence the
join $\langle W',W''\rangle$ in the geometric lattice $\Pi_W$ has rank either
$\rk(W^0)+3$ or $\rk(W^0)+4$. Conjugating by a suitable element of $W$,
we can assume that $W^0=W_J$, $W'=W_K$, and $\langle W',W''\rangle=W_L$,
where $J\subset K\subset L\subseteq I$,
$|K\setminus J|=2$, $|L\setminus K|=1$ or $2$,
and $J$ and $K$ have the property that all connected components are even.

If $|L\setminus K|=1$ and $L=I$, then $W_K$ and $W''$ have corank $1$ in $W$,
and so must be maximal in $\Pi_W^{(2)}$. If $|L\setminus K|=1$ and $L\neq I$,
let $i\in I\setminus L$ be any vertex adjacent to the connected 
component of $L$ which contains
the single vertex of $L\setminus K$. Then $L\cup\{i\}$ has all connected
components even, and $W_{L\cup\{i\}}$ is the required cover.

If $|L\setminus K|=2$, we must show that $L$ itself has the property
that all connected components are even. Suppose for a contradiction
that it did not. Then it must have exactly two odd connected components,
$L_1$ and $L_2$, containing the two elements $l_1$ and $l_2$ respectively
of $L\setminus K$. Clearly the two elements
of $K\setminus J$ are contained in the same connected component of $K$,
so (swapping $L_1$ and $L_2$ if necessary) we can assume that
$L_1\cap J=L_1\cap K=L_1\setminus\{l_1\}$. Now $W''$ is a parabolic
subgroup of $W_L$, hence a product of parabolic subgroups of the irreducible
components of $W_L$. The component of $W''$ inside $W_{L_1}$
must contain $W_{L_1\cap J}$ and cannot have odd rank, so it equals
$W_{L_1\cap J}$, which is also the component of $W_K$ inside $W_{L_1}$.
This contradicts the assumption that $\langle W_K,W''\rangle=W_L$.

Finally we prove (3). Since the product of two Cohen--Macaulay posets with
minimum elements is again Cohen--Macaulay (see \cite[Theorem 5.1.5 and
Exercise 5.1.6]{wachs}),
we may assume that $W$ is irreducible. If $\rk(W)$ is even, then $\Pi_W^{(2)}$
is a pure bounded finite poset which is totally semimodular, because every lower
order ideal $[\hat{0},W']$ in $\Pi_W^{(2)}$ is isomorphic to the product
of posets $\Pi_{W''}^{(2)}$ where $W''$ is irreducible of even rank.
So $\Pi_W^{(2)}$ satisfies the Bj\"orner--Wachs criterion for being 
Cohen--Macaulay (see \cite[Theorems 4.2.2 and 4.2.3]{wachs}). If $\rk(W)$
is odd, the same criterion shows that $\Pi_W^{(2)}\cup\{\hat{1}\}$
is Cohen--Macaulay, which implies that $\Pi_W^{(2)}$ is.
\epf

Thus in the Coxeter case, Theorem \ref{rainsthm} takes the
pleasing form:
\bth \label{rainscoxeterthm}
For all $i$, there is an isomorphism of $\Q$-vector spaces
\[ H^i(\overline{Y}_W(\R);\Q)\cong 
\bigoplus_{\substack{W'\in\Pi_W^{(2)}\\\rk(W')=2i}}
H_{i}(\hat{0},W')\otimes\orient((V^{W'})^\perp) \]
which is equivariant for $N_{GL(V)}(W)$.
\eth
\bpf
By Cohen--Macaulayness, the only poset homologies which survive on the
right-hand side of Theorem \ref{rainsthm} are those where
$\dim A-i=\rk_{\Pi_\CalG^{(2)}}(A)$, or in other words 
$i=\rk(W')-\rk_{\Pi_W^{(2)}}(W')
=\frac{\rk(W')}{2}$.
\epf
Note that if it were not for the twisting by $\orient((V^{W'})^\perp)$, the right-hand side of Theorem \ref{rainscoxeterthm} would be the
$i$th Whitney homology $\WH_i(\Pi_W^{(2)})$.

The $i=1$ special case of Theorem \ref{rainscoxeterthm} is simpler,
since $H_1(\hat{0},W')$ is canonically isomorphic to $\Q$ when $W'$ is an atom
of the poset. Hence we have an isomorphism of representations 
of $N_{GL(V)}(W)$:
\beq \label{h1eqn}
H^1(\overline{Y}_W(\R);\Q)\cong
\bigoplus_{W'}\Ind_{N_{GL(V)}(W')\,\cap\, N_{GL(V)}(W)}^{N_{GL(V)}(W)}
(\orient((V^{W'})^\perp)),
\eeq
where we sum over a set of representatives of the 
$N_{GL(V)}(W)$-conjugacy classes of irreducible rank-$2$
parabolic subgroups of $W$.
In the classical types this 
implies the following isomorphisms:
\beq \label{classh1eqn}
\begin{split}
H^1(\overline{Y}_{A_n}(\R);\Q)&\cong
\Ind_{S_{n-2}\times S_3}^{S_{n+1}}(1\boxtimes\varepsilon),\\
H^1(\overline{Y}_{B_n}(\R);\Q)&\cong
\Ind_{W(B_{n-2})\times W(B_2)}^{W(B_n)}(1\boxtimes\varepsilon)\\
&\qquad\oplus\Ind_{W(B_{n-3})\times \{\pm 1\}\times S_3}^{W(B_n)}(1\boxtimes
1\boxtimes\varepsilon),\\
H^1(\overline{Y}_{D_n}(\R);\Q)&\cong
\Ind_{W(B_{n-3})\times \{\pm 1\}\times S_3}^{W(B_n)}(1\boxtimes
1\boxtimes\varepsilon).
\end{split}
\eeq
To describe the higher cohomology groups, we need further information
about the homology of the poset $\Pi_W^{(2)}$; this will be the focus
of the next two sections.

To save repetition, we note a close connection between the
Whitney homology of $\Pi_{B_n}^{(2)}$ and $\Pi_{D_n}^{(2)}$,
which should be compared with the relationship between
the Whitney homology of the lattices $\Pi_{B_n}$ and $\Pi_{D_n}$
(see \cite[Theorem 5.4]{hendersondowling}).
Use the standard interpretations $B_0=$ empty arrangement, $B_1=A_1$,
$D_2=A_1\times A_1$, $D_3=A_3$.
\bpr \label{typedprop}
For all $i$ and $n\geq 2$, we 
have an isomorphism
\[ \WH_i(\Pi_{B_n}^{(2)})\cong \WH_i(\Pi_{D_n}^{(2)})
\oplus
\Ind_{W(B_{n-2})\times W(B_2)}^{W(B_n)}(\WH_{i-1}(\Pi_{B_{n-2}}^{(2)})
\boxtimes 1) \]
of representations of $W(B_n)$. 
\epr
\bpf
The general result \cite[Theorem 3.2]{rainshomology}
provides a canonical isomorphism $H_i(0,A)\cong H_i^f(A)$ for all
$A\in\Pi_{\CalG}^{(2)}$, where $H_i^f(A)$ denotes the homology
of the chain complex $C_\bullet^f(A)$ of $2$-divisible forests with root $A$.
With the usual definition of the hyperplane arrangement
$\CalA_{B_n}$, the building set $\CalG_{B_n}$ consists of the 
following subspaces,
written using a basis $x_1,\cdots,x_n$ of $V^*$:
\bes
\begin{split} 
&\R^I=\{\sum_{i\in I}a_i x_i\},\text{ for $I\subseteq\{1,\cdots,n\}$, $|I|\geq 1$,
and}\\
&\R_\varphi^L=
\{\sum_{\ell\in L}a_\ell x_\ell\,|\,
\sum_{\ell\in L}\varphi(\ell)a_\ell=0\},\\
&\qquad\text{ for $L\subseteq\{1,\cdots,n\}$, $|L|\geq 2$, 
$\varphi:L\to\{\pm 1\}$.}
\end{split}
\ees
The building set $\CalG_{D_n}$ can be identified with the subset
of $\CalG_{B_n}$ consisting of all but the subspaces $\R^I$ where
$|I|\leq 2$. From this it is clear that
any $2$-divisible forest for $\CalG_{D_n}$ is also a $2$-divisible
forest for $\CalG_{B_n}$; the only $2$-divisible forests
for $\CalG_{B_n}$ which are not of this kind are those which contain
a subspace $\R^{\{i,j\}}$, necessarily unique. Hence we have a short
exact sequence of chain complexes:
\[ 0 \to C_\bullet^{f,D_n}(A) \to C_\bullet^{f,B_n}(A) \to
\bigoplus_{\R^{\{i,j\}}\subseteq A} C_{\bullet-1}^{f,B_{n-2}}(A/\R^{\{i,j\}})
\to 0 \]
for all $A\in\Pi_{\CalG_{B_n}}^{(2)}$ (interpret 
$C_\bullet^{f,D_n}(A)$ as $0$
if $A\not\in\Pi_{\CalG_{D_n}}^{(2)}$). The corresponding long exact sequence
of homology groups reduces to a short exact sequence, because of the
Cohen-Macaulay property proved in Theorem \ref{posetthm}. The result follows.
\epf
Note that the direct sum decomposition in Proposition \ref{typedprop} is not
canonical, since the short exact sequence has no canonical splitting.
\bcr \label{typedcor}
For all $i$ and $n\geq 2$, we have an isomorphism
\bes
\begin{split} 
H^i(\YBn(\R);\Q)&\cong H^i(\YDn(\R);\Q)\\
&\quad\oplus
\Ind_{W(B_{n-2})\times W(B_2)}^{W(B_n)}
(H^{i-1}(\overline{Y}_{B_{n-2}}(\R);\Q)\boxtimes
\varepsilon)
\end{split}
\ees
of representations of $W(B_n)$. 
\ecr
\bpf
This follows from Theorem \ref{rainscoxeterthm}, by incorporating
the factors $\orient((V^{W'})^\perp)$ into the argument of the
Proposition.
\epf
To be more concrete, for $n\geq 3$
the inclusion of $\CalG_{D_n}$ in
$\CalG_{B_n}$ induces a surjective birational map $\YBn\to\YDn$.
In fact, $\YBn$ is the blow-up of $\YDn$ along a collection of
disjoint subvarieties, each isomorphic to 
$\overline{Y}(B_{n-2}(\R))$, indexed by the
cosets $W(B_n)/(W(B_{n-2})\times W(B_2))$. The induced map
in cohomology $H^i(\YDn(\R);\Q)\to H^i(\YBn(\R);\Q)$ is the
$W(B_n)$-equivariant injection seen in Corollary \ref{typedcor}.
\section{Betti numbers and Euler characteristic}
In this section we neglect the action of $N_{GL(V)}(W)$, and seek
to compute merely the Betti numbers and Euler characteristic of
$\overline{Y}_W(\R)$. We encode the Betti numbers in the Poincar\'e
polynomial
\[ P(\overline{Y}_W(\R),t)=\sum_{i\geq 0} 
\dim H^i(\overline{Y}_W(\R);\Q)\, (-t)^i, \]
so that the Euler characteristic
$\chi(\overline{Y}_W(\R))$ equals $P(\overline{Y}_W(\R),1)$. As a
consequence of Theorem \ref{rainscoxeterthm} we have:
\bpr
The Poincar\'e polynomial $P(\overline{Y}_W(\R),t)$
equals the characteristic polynomial
\[ \chi(\Pi_W^{(2)},t)=
\sum_{W'\in\Pi_W^{(2)}}\mu_{\Pi_W^{(2)}}(\hat{0},W')\,t^{\rk(W')/2}, \]
where $\mu_{\Pi_W^{(2)}}$ is the M\"obius function.
In particular, the
Euler characteristic $\chi(\overline{Y}_W(\R))$ equals
$\chi(\Pi_W^{(2)},1)$.
\epr
\bpf
Take dimension of both sides of Theorem \ref{rainscoxeterthm}, and
use the fact that $\Pi_W^{(2)}$ is Cohen--Macaulay.
\epf
Computing $P(\overline{Y}_W(\R),t)$ is thus reduced
to the purely combinatorial problem of
finding $\chi(\Pi_W^{(2)},t)$.
Clearly this is $1$ if $n=1$ and $1-t$ if
$n=2$. Moreover,
\beq
\chi(\Pi_W^{(2)},1)=0\text{ when $n$ is even,}
\eeq
because the poset then has
a maximum (correspondingly, the Euler characteristic of
$\overline{Y}_W(\R)$ is zero because it is an odd-dimensional compact
manifold).
This, and our knowledge that
the coefficient of $-t$ is the number of irreducible rank-$2$
parabolic subgroups, is enough information to handle the
exceptional types of rank $\leq 4$:
\beq
\begin{split}
\chi(\Pi_{H_3}^{(2)},t)&=1-16t,\\
\chi(\Pi_{F_4}^{(2)},t)&=1-50t+49t^2,\\
\chi(\Pi_{H_4}^{(2)},t)&=1-272t+271t^2.
\end{split}
\eeq
Using the technique of Sundaram (see \cite[Theorem 4.4.1]{wachs})
to simplify the calculations, one can compute:
\beq
\begin{split}
\chi(\Pi_{E_6}^{(2)},t)&=1-120t+2739t^2-2620t^3,\\
\chi(\Pi_{E_7}^{(2)},t)&=1-336t+26229t^2-230014t^3,\\
\chi(\Pi_{E_8}^{(2)},t)&=1-1120t+332178t^2-18066280t^3+17735221t^4.
\end{split}
\eeq
For this we used the Magma computational algebra package.

It only remains to treat the classical types.
In \cite{chr} (see also \cite[Corollary 3.16]{etingofetal}) it is shown that
\beq
\chi(\Pi_{A_n}^{(2)},t)=\prod_{k=1}^{\lfloor\frac{n}{2}\rfloor}
(1-(n+1-2k)^2 t),
\eeq
which implies that
\beq \label{typeaeulercharacteristiceqn}
\chi(\Pi_{A_n}^{(2)},1)=
\frac{(-1)^{\frac{n-1}{2}}(n+1)!(n-1)!}
{2^n (\frac{n+1}{2})!(\frac{n-1}{2})!},\text{ when $n$ is odd.}
\eeq
This formula for $\chi(\overline{Y}_{A_n}(\R))$ was proved topologically
by Devadoss (\cite[Theorem 3.2.3]{devadoss}).

A generating-function formula for the characteristic polynomial in type $B$
was effectively found in \cite{hendersondowling}:
\bpr \label{typebnoneqprop}
In $\Q[t][\![x]\!]$ we have
\bes
\begin{split} 
1 + \frac{x}{2} + \sum_{n\geq 2}\,
& \chi(\Pi_{B_n}^{(2)},t)\frac{x^n}{2^n n!}\\
&=\sech(\frac{1}{2}\arcsinh(t^{1/2}x))
\exp(\frac{t^{-1/2}}{2}\arcsinh(t^{1/2}x)).
\end{split}
\ees
\epr
\bpf
This is the $G=\{\pm 1\}$ case of \cite[(5.10)]{hendersondowling},
with the first term on the right-hand side omitted since it 
corresponds to the maximum elements which are not included
in $\Pi_{B_n}^{(2)}$.
\epf
Using the identity 
$\tanh(\frac{1}{2}\arcsinh(x))
=\frac{(1+x^2)^{1/2}-1}{x}$, we deduce that
\beq
\chi(\Pi_{B_n}^{(2)},1)=
\frac{(-1)^{\frac{n-1}{2}}n!(n-1)!}
{(\frac{n+1}{2})!(\frac{n-1}{2})!},\text{ when $n$ is odd.}
\eeq
The fact that $\chi(\overline{Y}_{B_n}(\R))=\frac{2^n}{n+1}
\chi(\overline{Y}_{A_n}(\R))$ was proved topologically by Gaiffi
(see \cite[Theorem 6.9]{gaiffi}).
\bcr
For $n\geq 2$, $\chi(\Pi_{B_n}^{(2)},t)$ equals
\[ \sum_{m=0}^{\lfloor\frac{n}{2}\rfloor}
\binom{n}{2m}\frac{(4m)!}{2^{2m}(2m+1)!}\,(-t)^m
\negthickspace\negthickspace
\prod_{\substack{1\leq a\leq n-2m-2\\a\equiv n\text{ mod }2}}
\negthickspace\negthickspace
(1-4a^2 t). \]
\ecr
\bpf
By a standard identity,
\beq
\sech(\frac{1}{2}\arcsinh(t^{1/2}x))=\sum_{m\geq 0}\frac{(4m)!}
{2^{2m}(2m+1)!}\,(-t)^m \frac{x^{2m}}{2^{2m}(2m)!},
\eeq 
and by an argument entirely analogous to
\cite[Proof of Corollary 3.16]{etingofetal},
\beq
\exp(\frac{t^{-1/2}}{2}\arcsinh(t^{1/2}x))
=\sum_{m\geq 0}\prod_{\substack{1\leq a\leq m-2\\a\equiv m\text{ mod }2}}
(1-4a^2 t)\frac{x^m}{2^m m!}.
\eeq
So the result follows from Proposition \ref{typebnoneqprop}.
\epf
In contrast to type $A$, $\chi(\Pi_{B_n}^{(2)},t)$
does not generally factorize into linear factors
with integer coefficients, e.g.\
$\chi(\Pi_{B_5}^{(2)},t)=1-50t+289t^2$.

Finally, Corollary \ref{typedcor} shows that
\beq
\chi(\Pi_{D_n}^{(2)},t)=
\chi(\Pi_{B_n}^{(2)},t)+\binom{n}{2}t\chi(\Pi_{B_{n-2}}^{(2)},t),
\eeq
which implies that
\beq
\chi(\Pi_{D_n}^{(2)},1)=
\frac{(-1)^{\frac{n-1}{2}}(n-1)(7n-17)\, n!(n-3)!}
{8(\frac{n+1}{2})!(\frac{n-1}{2})!},
\text{ when $n\geq 3$ is odd.}
\eeq
The fact that $\chi(\overline{Y}_{D_n}(\R))=\frac{2^{n-3}(7n-17)}{(n+1)(n-2)}
\chi(\overline{Y}_{A_n}(\R))$ was proved topologically in
\cite[Theorem 6.14]{armstrongetal}.
\section{The $W(B_n)$-action on the cohomology of $\YBn(\R)$}
We now return to the equivariant setting, and deduce from
Theorem \ref{rainscoxeterthm} a
plethystic generating-function formula
for the characters of $W(B_n)$ acting
on the rational cohomology groups of $\YBn(\R)$, analogous
to the type-$A$ formula given in \cite{rainstypea}. (Because of
Corollary \ref{typedcor}, any information about $\YBn(\R)$
has immediate implications for $\YDn(\R)$.)

Let $\Lambda_A$ denote
the usual $\N$-graded ring of symmetric functions, with scalars extended to
$\Q$; thus $\Lambda_A=\Q[p_i\,|\,i\geq 1]$ where $p_i$ is the
power sum function of degree $i$. The character of a representation
$M$ of $S_n$ over $\Q$ is encapsulated in its Frobenius characteristic
$\ch_{S_n}(M)$, which is a homogeneous element of $\Lambda_A$ of degree $n$.
Let ${}^\sign$ denote the ring involution on $\Lambda_A$ which
corresponds
to tensoring representations with the sign character $\varepsilon$; 
that is, the one
satisfying $p_i^\sign=(-1)^{i-1}p_i$.
We also need the associative operation $\circ$ on $\Lambda_A$ called plethysm
(see \cite[Chapter I, Section 8]{macdonald}), which satisfies
$p_i\circ p_j=p_{ij}$; in particular, $p_1$ is the plethystic identity.
Note that if $g\in\Lambda_A$ has only odd-degree terms,
$(f\circ g)^\sign=f^\sign\circ g^\sign$ for all $f\in\Lambda_A$.

Following \cite[Chapter I, Appendix B]{macdonald}, we introduce a type-$B$
analogue $\Lambda_B:=\Q[x_i,y_i\,|\,i\geq 1]$,
which is $\N$-graded by setting $\deg(x_i)=\deg(y_i)=i$.
(In the notation of [loc.\ cit.], $W(B_n)$ is $\{\pm 1\}\wr S_n$,
$x_i$ is $p_i(\{1\})$ and $y_i$ is $p_i(\{-1\})$.)
The Frobenius characteristic of a representation $M$ of $W(B_n)$ over $\Q$ is
\beq
\ch_{W(B_n)}(M):=\frac{1}{2^n n!}\sum_{w\in W(B_n)}
\tr(w,M)\,\prod_{i\geq 1} x_i^{a_i(w)}y_i^{b_i(w)},
\eeq
where $a_i(w)$ (respectively, $b_i(w)$)
denotes the number of cycles in $w$ of length $i$ whose cycle-product is 
$1$ (respectively, $-1$). This
is a homogeneous element of $\Lambda_B$ of degree $n$.
We define a ring involution ${}^\sign$ on $\Lambda_B$ by the rules
$x_i^\sign=(-1)^{i-1}x_i$, $y_i^{\sign}=(-1)^i y_i$; again this
corresponds to tensoring representations with the sign character 
$\varepsilon$ of $W(B_n)$.
We also have a right plethystic action
$\circ:\Lambda_B\times\Lambda_A\to\Lambda_B$, which satisfies 
\beq
x_i\circ p_j=x_{ij},\ y_i\circ p_j=\left\{\begin{array}{cl}
x_{ij},\text{ if $j$ is even,}\\
y_{ij},\text{ if $j$ is odd.}\end{array}\right.
\eeq
If $g\in\Lambda_A$ has only odd-degree terms,
$(f\circ g)^\sign=f^\sign\circ g^\sign$ for all $f\in\Lambda_B$.
See \cite[Section 5]{hendersonwreath} for some other properties of this
action.

Complete $\Lambda_A$ and $\Lambda_B$ to the
corresponding formal power series rings
$\Aa_A=\Q[\![p_i]\!]$ and $\Aa_B=\Q[\![x_i,y_i]\!]$.
We can extend plethysm and the
plethystic action to these in the obvious way, provided that
the right-hand
input has zero constant term. We define
\bes
\begin{split}
\Exp&:=\sum_{n\geq 0}\ch_{S_n}(\bbone)
=\exp(\sum_{m\geq 1}\frac{p_m}{m})\ \in\Aa_A,\\
\Exp_B&:=\sum_{n\geq 0}\ch_{W(B_n)}(\bbone)
=\exp(\sum_{m\geq 1}\frac{x_m+y_m}{2m})\ \in\Aa_B.
\end{split}
\ees
We also define $\Cosh$ and $\Sinh$ (respectively, $\Cosh_B$ and $\Sinh_B$) 
to be the sums of the even-degree and odd-degree terms of $\Exp$
(respectively, $\Exp_B$).
Since $\Cosh_B$ has constant term $1$, it has a multiplicative inverse
$\Sech_B\in\Aa_B$, which also has only even-degree terms. 
Since $\Sinh$ has zero constant term, it has a
plethystic inverse $\Arcsinh\in\Aa_A$, which has only odd-degree terms.

To keep track of the grading by cohomological degree,
we must work in $\Aa_A\otimes_\Q \Q[t,t^{-1}]$ and
$\Aa_B\otimes_\Q \Q[t,t^{-1}]$ where $t$ is another indeterminate.
We extend the plethysm and plethystic action to these rings by the rules
$p_i\circ t=x_i\circ t=y_i\circ t=t^i$. Actually the formulae can be
stated more conveniently if we use $\Q[t^{1/2},t^{-1/2}]$ rather than
just $\Q[t,t^{-1}]$, although all exponents of $t$ in the final
result are integers.

In this notation, the generating function for the character
of $S_n$ acting on the rational cohomology groups of 
$\overline{Y}_{A_{n-1}}(\R)$ is as follows.
\bth \cite[Theorem 3.5]{rainstypea} \label{rainstypeathm}
In $\Aa_A\otimes_\Q \Q[t,t^{-1}]$ we have
\bes
\begin{split} 
1+p_1+\sum_{n\geq 2}\sum_{i\geq 0}
\ch_{S_n}&(H^i(\overline{Y}_{A_{n-1}}(\R);\Q))\,(-t)^i\\
&= \Exp\circ t^{-1/2}\Arcsinh^\sign\circ t^{1/2}p_1.
\end{split}
\ees
\eth
\noindent
(Here $t^{-1/2}\Arcsinh^\sign\circ t^{1/2}p_1$ is an element
of $\Aa_A\otimes_\Q \Q[t,t^{-1}]$ because all terms of $\Arcsinh^\sign$
have odd degree; in fact, this is just $\Arcsinh^\sign$ with the
degree-$(2m+1)$ term multiplied by $t^m$.)

The type-$B$ analogue is:
\bth \label{mainthm}
In $\Aa_B\otimes_\Q \Q[t,t^{-1}]$ we have
\bes
\begin{split} 
1 + &\frac{x_1+y_1}{2} + \sum_{n\geq 2}\sum_{i\geq 0}
\ch_{W(B_n)}(H^i(\YBn(\R);\Q))\,(-t)^i\\
&=(\Sech_B^\sign\circ\Arcsinh^\sign\circ t^{1/2}p_1)
(\Exp_B\circ t^{-1/2}\Arcsinh^\sign\circ t^{1/2}p_1).
\end{split}
\ees
\eth
\bpf
By Theorem \ref{rainscoxeterthm}, we can identify
$H^\bullet(\YBn(\R);\Q)$ with the Whitney homology of the poset
$\Pi_{B_n}^{(2)}$, twisted by the $\orient((V^{W'})^\perp)$ factors.
Thus the proof proceeds almost identically to the proof for the untwisted
Whitney homology of $Q_n^{1\, \mathrm{mod}\, 2}(\{\pm 1\})$ given
in \cite[Section 5]{hendersondowling}. The key point there is that
the lower order ideal $[\hat{0},(J,\pi)]$ is isomorphic to the
poset product 
$Q^{1\text{ mod }2}(J)\times\prod_{K\in\pi'}\Pi^{1\text{ mod }2}(K)$, 
where the poset-valued
functors $Q^{1\text{ mod }2}$ and $\Pi^{1\text{ mod }2}$ are defined
as in \cite[Section 4]{hendersondowling}, and $\pi'$ denotes
a set of representatives for the $\{\pm 1\}$-orbits on $\pi$. So by
the K\"unneth formula \cite[second statement of Theorem 5.1.5]{wachs},
\beq 
H_{\rk(J,\pi)}(\hat{0},(J,\pi))\cong\tH_{Q^{1\text{ mod }2}}(J)
\otimes\bigotimes_{K\in\pi'}\tH_{\Pi^{1\text{ mod }2}}(K),
\eeq
where the graded $(\{\pm 1\}\wr \Ss)$-module $\tH_{Q^{1\text{ mod }2}}$ 
and the graded
$\Ss$-module $\tH_{\Pi^{1\text{ mod }2}}$ are defined as in 
\cite[Section 4]{hendersondowling}. 
The new ingredient here is $\orient((V^{W'})^\perp)$; it is easy to
see that if $W'\in\Pi_{B_n}^{(2)}$ corresponds to 
$(J,\pi)\in Q_n^{1\, \mathrm{mod}\, 2}(\{\pm 1\})$, then
\beq 
\orient((V^{W'})^\perp)\cong\varepsilon_{\{\pm 1\}}(J)\otimes
\bigotimes_{K\in\pi'}\varepsilon(K),
\eeq
where $\varepsilon$ is the $\Ss$-module such that
$\varepsilon(K)$ is a one-dimensional vector space on which the
permutations of $K$ act via the sign character, and
$\varepsilon_{\{\pm 1\}}$ is the analogous $(\{\pm 1\}\wr \Ss)$-module.
So we can translate Theorem \ref{rainscoxeterthm} into the 
following equation in $\Aa_B\otimes_\Q \Q[t,t^{-1}]$:
\beq \label{intermediateeqn}
\begin{split} 
1 + & \frac{x_1+y_1}{2} + \sum_{n\geq 2}\sum_{i\geq 0}
\ch_{W(B_n)}(H^i(\YBn(\R);\Q))\,(-t)^i\\
&\quad=\ch_t[(\tH_{Q^{1\text{ mod }2}}^{0\text{ mod }2}\otimes
\varepsilon_{\{\pm 1\}})\cdot(\bbone_{\{\pm 1\}}\circ
(\tH_{\Pi^{1\text{ mod }2}}^{1\text{ mod }2}\otimes\varepsilon))].
\end{split}
\eeq
As well as the notation $\ch_t$, $\cdot$, $\circ$ introduced in 
\cite{hendersondowling}, we have here used the obvious definition
of tensor-multiplying a graded $\Ss$-module (respectively,
a graded $(\{\pm 1\}\wr \Ss)$-module)
by the $\Ss$-module 
$\varepsilon$ (respectively, the $(\{\pm 1\}\wr \Ss)$-module
$\varepsilon_{\{\pm 1\}}$).
Now we apply \cite[Theorem 5.1]{hendersondowling} to 
the right-hand side of \eqref{intermediateeqn}, which converts
it into
\[ \ch_t(\tH_{Q^{1\text{ mod }2}}^{0\text{ mod }2})^\sign\cdot
(\Exp_B\circ\ch_t(\tH_{\Pi^{1\text{ mod }2}}^{1\text{ mod }2})^\sign). \]
The $d=2$ case of \cite[(5.9)]{hendersondowling}, a rephrasing of
\cite[Theorem 4.7]{chr}, says that
\[ \ch_t(\tH_{\Pi^{1\text{ mod }2}}^{1\text{ mod }2})=t^{-1/2}\Arcsinh\circ
t^{1/2}p_1, \]
and the $d=2$ case of \cite[Theorem 2.6]{hendersondowling}
says that
\[ \ch_t(\tH_{Q^{1\text{ mod }2}}^{0\text{ mod }2})=
\Sech_B\circ\Arcsinh\circ t^{1/2}p_1. \]
Theorem \ref{mainthm} follows.
\epf
The non-equivariant Proposition \ref{typebnoneqprop} can be obtained
from this Theorem by setting $x_1\to x$ and all other $x_i,y_i\to 0$.

The two factors in the right-hand side of Theorem \ref{mainthm}
have neatly separable meanings.
View $W(B_n)$ as $\Gamma_n\rtimes S_n$ where
$\Gamma_n=\{\pm 1\}^n$, and 
note that for any representation $M$ of $W(B_n)$,
$S_n$ preserves the $\Gamma_n$-invariant subspace $M^{\Gamma_n}$ and also the
subspace $M^{\Gamma_n,\rmprod}$ on which $\Gamma_n$ acts via its character
$\rmprod:(\varepsilon_1,\cdots,\varepsilon_n)\mapsto\varepsilon_1\cdots
\varepsilon_n$.
\bpr \label{normalsubgpprop}
In $\Aa_A\otimes_\Q \Q[t,t^{-1}]$ we have
\bes
\begin{split} 
1+p_1+\sum_{n\geq 2}\sum_{i\geq 0}
\ch_{S_n}&(H^i(\YBn(\R);\Q)^{\Gamma_n})\,(-t)^i\\
&\qquad=\Exp\circ t^{-1/2}\Arcsinh^\sign\circ t^{1/2}p_1,\\
1+\sum_{n\geq 2}\sum_{i\geq 0}
\ch_{S_n}&(H^i(\YBn(\R);\Q)^{\Gamma_n,\rmprod})\,(-t)^i\\
&\qquad=\Sech^\sign\circ\Arcsinh^\sign\circ t^{1/2}p_1.
\end{split}
\ees
\epr
\bpf
It is easy to see that for any representation $M$ of $W(B_n)$,
\beq 
\ch_{S_n}(M^{\Gamma_n})=\ch_{W(B_n)}(M)|_{x_i,y_i\to p_i}. 
\eeq
So to prove the first equation we apply the substitution $x_i,y_i\to p_i$
to the right-hand side of Theorem \ref{mainthm}, noting that
it respects the plethystic action. The first factor becomes
$(\Sech_B^\sign|_{x_i,y_i\to p_i})\circ\Arcsinh^\sign\circ t^{1/2}p_1$,
which equals $1$ because
\beq \label{expsigneqn}
\Exp_B^\sign|_{x_i,y_i\to p_i}
=\Exp_B|_{x_i\to(-1)^{i-1}p_i, 
y_i\to(-1)^i p_i}
=1,
\eeq
by definition of $\Exp_B$.
The second factor becomes the right-hand side of the stated result, because
\beq \label{expeqn} 
\Exp_B|_{x_i,y_i\to p_i}
=\Exp.
\eeq
Similarly, for any representation $M$ of $W(B_n)$,
\beq
\ch_{S_n}(M^{\Gamma_n,\rmprod})=\ch_{W(B_n)}(M)|_{x_i\to p_i, y_i\to -p_i}.
\eeq
So to prove the second equation we apply this substitution
to the right-hand side of Theorem \ref{mainthm}. The first factor becomes
\bes
\begin{split} 
(\Sech_B\circ\Arcsinh\circ t^{1/2}p_1)^\sign&|_{x_i\to p_i, y_i\to -p_i}\\
&=(\Sech_B\circ\Arcsinh\circ t^{1/2}p_1)|_{x_i,y_i\to(-1)^{i-1}p_i}\\
&=((\Sech_B\circ\Arcsinh\circ t^{1/2}p_1)|_{x_i,y_i\to p_i})^\sign\\
&=(\Sech\circ\Arcsinh\circ t^{1/2}p_1)^\sign,
\end{split}
\ees
which is the right-hand side of the stated result (in the last step
we have used \eqref{expeqn} again). The second factor
becomes
\bes
\begin{split}
(\Exp_B\circ t^{-1/2}\Arcsinh^\sign&\circ t^{1/2}p_1)|_{x_i\to p_i, 
y_i\to -p_i}\\
&=(\Exp_B^\sign\circ t^{-1/2}\Arcsinh\circ t^{1/2}p_1)|_{x_i,y_i\to
(-1)^{i-1} p_i}\\
&=((\Exp_B^\sign\circ t^{-1/2}\Arcsinh\circ t^{1/2}p_1)|_{x_i,
y_i\to p_i})^\sign\\
&=1,
\end{split}
\ees
where in the last step we have used \eqref{expsigneqn} again.
\epf
Note that the first part of this Proposition, together with
Theorem \ref{rainstypeathm}, implies 
an isomorphism of graded $S_n$-representations
\beq \label{isomeqn}
H^\bullet(\YBn(\R);\Q)^{\Gamma_n}\cong 
H^\bullet(\overline{Y}_{A_{n-1}}(\R);\Q).
\eeq
It is natural to suspect that this isomorphism is induced by
a homotopy equivalence between 
$\Gamma_n\!\setminus\!\YBn(\R)$ (the topological quotient) and 
$\overline{Y}_{A_{n-1}}(\R)$. By comparison, the formula analogous
to Theorem \ref{mainthm} for the complex manifold $\YBn(\C)$,
which is the $r=2$ case of \cite[Theorem 3.4]{hendersonwreath},
implies an isomorphism of graded $S_n$-representations
\beq \label{complexisomeqn}
H^\bullet(\YBn(\C);\Q)^{\Gamma_n}\cong 
H^\bullet(\overline{Y}_{A_{n}}(\C);\Q),
\eeq
where the right-hand side is restricted from $S_{n+1}$ to $S_n$.

Specializing further, we can find the occurrences
of the four one-dimensional characters of $W(B_n)$ in
$H^\bullet(\YBn(\R);\Q)$. For this we need a preparatory result:
\blm \label{arcsinhlemma}
In $\Q[\![x]\!]$,
$\Arcsinh|_{p_i\to x^i}=x-x^3$, $\Arcsinh^\sign|_{p_i\to x^i}=x$.
\elm
\bpf
The first claim is proved in \cite[p. 306]{chr}, and the proof of the
second is similar.
\epf
Let $\rmprod$ denote the one-dimensional character
of $W(B_n)$ whose restriction to $\Gamma_n$ is $\rmprod$ as above, and
whose restriction to $S_n$ is trivial.
\bpr \label{onedimcharprop}
For $n\geq 2$, the multiplicities of the one-dimensional characters 
of $W(B_n)$ in $H^i(\YBn(\R);\Q)$ are:
\bes
\begin{split}
\lng H^i(\YBn(\R);\Q), 1\rng_{W(B_n)}
&=\left\{\begin{array}{cl}
1,&\text{ if $i=0$,}\\
0,&\text{ otherwise,}
\end{array}\right.\\
\lng H^i(\YBn(\R);\Q), \varepsilon.\rmprod\rng_{W(B_n)}
&=\left\{\begin{array}{cl}
1,&\text{ if $n\not\equiv 2$ mod $3$, $i=\lfloor\frac{n}{3}\rfloor$,}\\
0,&\text{ otherwise,}
\end{array}\right.\\
\lng H^i(\YBn(\R);\Q), \rmprod\rng_{W(B_n)}
&=0,\\
\lng H^i(\YBn(\R);\Q), \varepsilon\rng_{W(B_n)}
&=\left\{\begin{array}{cl}
1,&\text{ if $n\equiv 0$ mod $2$, $i=\frac{n}{2}$,}\\
0,&\text{ otherwise.}
\end{array}\right.
\end{split}
\ees
In particular, the topological quotient
$W(B_n)\!\setminus\!\YBn(\R)$ has the rational cohomology of a point.
\epr
\bpf
For any representation $M$ of $S_n$, 
$\ch_{S_n}(M)|_{p_i\to x^i}=\lng M, 1\rng_{S_n}\, x^n$.
Now the first multiplicity in the statement is the
multiplicity of the trivial character of $S_n$ in
$H^i(\YBn(\R);\Q)^{\Gamma_n}$; according to the first part of
Proposition \ref{normalsubgpprop},
we find this by computing
\bes
\begin{split}
(\Exp\circ t^{-1/2}&\Arcsinh^\sign\circ t^{1/2}p_1)|_{p_i\to x^i}\\
&=\exp(\sum_{m\geq 1}\frac{1}{m}(p_m\circ t^{-1/2}\Arcsinh^\sign
\circ t^{1/2}p_1)|_{p_i\to x^{i}})\\
&=\exp(\sum_{m\geq 1}\frac{1}{m}(t^{-1/2}\Arcsinh^\sign\circ t^{1/2}p_1)|_{p_i\to x^{mi},\,t\to t^m})\\
&=\exp(\sum_{m\geq 1}\frac{1}{m}x^m)\\
&=1+x+x^2+x^3+\cdots,
\end{split}
\ees
where the third equality uses Lemma \ref{arcsinhlemma}.
Similarly, the second multiplicity is that of
the sign character of $S_n$ in
$H^i(\YBn(\R);\Q)^{\Gamma_n}$, which we find by computing
\bes
\begin{split}
(\Exp^\sign\circ t^{-1/2}&\Arcsinh\circ t^{1/2}p_1)|_{p_i\to x^i}\\
&=\exp(\sum_{m\geq 1}\frac{(-1)^{m-1}}{m}(p_m\circ t^{-1/2}\Arcsinh
\circ t^{1/2}p_1)|_{p_i\to x^{i}})\\
&=\exp(\sum_{m\geq 1}\frac{(-1)^{m-1}}{m}(t^{-1/2}\Arcsinh\circ 
t^{1/2}p_1)|_{p_i\to x^{mi},\,t\to t^m})\\
&=\exp(\sum_{m\geq 1}\frac{(-1)^{m-1}}{m}(x^m-t^m x^{3m}))\\
&=1+x-tx^3-tx^4+t^2 x^6+t^2 x^7-\cdots.
\end{split}
\ees
To treat the second part of Proposition \ref{normalsubgpprop} similarly,
note that 
\bes
\begin{split}
(\Exp^\sign\circ\Arcsinh^\sign\circ t^{1/2}p_1)|_{p_i\to x^i}
&=\exp(\sum_{m\geq 1}\frac{(-1)^{m-1}}{m}t^{m/2}x^m)\\
&=1+t^{1/2}x,
\end{split}
\ees
implying that
\[ (\Sech^\sign\circ\Arcsinh^\sign\circ t^{1/2}p_1)|_{p_i\to x^i}=1, \]
and
\bes
\begin{split}
(\Exp\circ\Arcsinh\circ t^{1/2}p_1)|_{p_i\to x^i}
&=\exp(\sum_{m\geq 1}\frac{1}{m}(t^{m/2}x^m-t^{3m/2} x^{3m}))\\
&=1+t^{1/2}x+t x^2,
\end{split}
\ees
implying that
\[ (\Sech\circ\Arcsinh\circ t^{1/2}p_1)|_{p_i\to x^i}=
(1+t x^2)^{-1}=1-t x^2+t^2 x^4-t^3 x^6+\cdots. \]
The third and fourth multiplicity statements follow.
\epf
It is easy to deduce the multiplicities of one-dimensional characters
of $W(B_n)$ and $W(D_n)$ in the cohomology groups of $\YDn(\R)$:
\bpr
For $n\geq 2$, we have
\bes
\begin{split}
\lng H^i(\YDn(\R);\Q), 1\rng_{W(B_n)}
&=\left\{\begin{array}{cl}
1,&\text{ if $i=0$,}\\
0,&\text{ otherwise,}
\end{array}\right.\\
\lng H^i(\YDn(\R);\Q), \varepsilon.\rmprod\rng_{W(B_n)}
&=\left\{\begin{array}{cl}
1,&\text{ if $n\not\equiv 2$ mod $3$, $i=\lfloor\frac{n}{3}\rfloor$,}\\
0,&\text{ otherwise,}
\end{array}\right.\\
\lng H^i(\YDn(\R);\Q), \rmprod\rng_{W(B_n)}
&=0,\\
\lng H^i(\YDn(\R);\Q), \varepsilon\rng_{W(B_n)}
&=0,\\
\lng H^i(\YDn(\R);\Q), 1\rng_{W(D_n)}
&=\left\{\begin{array}{cl}
1,&\text{ if $i=0$,}\\
0,&\text{ otherwise,}
\end{array}\right.\\
\lng H^i(\YDn(\R);\Q), \varepsilon\rng_{W(D_n)}
&=\left\{\begin{array}{cl}
1,&\text{ if $n\not\equiv 2$ mod $3$, $i=\lfloor\frac{n}{3}\rfloor$,}\\
0,&\text{ otherwise.}
\end{array}\right.
\end{split}
\ees
In particular, the topological quotient
$W(D_n)\!\setminus\!\YDn(\R)$
has the rational cohomology of a point.
\epr
\bpf
To obtain the first four statements, combine Proposition \ref{onedimcharprop}
and Corollary \ref{typedcor}. The other two statements follow because
$\Ind_{W(D_n)}^{W(B_n)}(1)=1+\rmprod$ and
$\Ind_{W(D_n)}^{W(B_n)}(\varepsilon)=\varepsilon+\varepsilon.\rmprod$.
\epf

Another consequence of Theorem \ref{mainthm} is a generating function
for the Euler character. First recall that by setting $t\to 1$ in
Theorem \ref{rainstypeathm} we obtain the following equation in $\Aa_A$:
\beq \label{typeaeulerchareqn}
1+p_1+\sum_{n\geq 2}\sum_{i\geq 0}(-1)^i
\ch_{S_n}(H^i(\overline{Y}_{A_{n-1}}(\R);\Q))
= \Exp\circ\Arcsinh^\sign.
\eeq
As observed in \cite{rainstypea}, the right-hand side is a power series
in the indeterminates $p_{2^l}$, $l\geq 0$, only. The type-$B$ analogue is:
\bpr \label{eulercharprop}
In $\Aa_B$ we have
\bes
\begin{split} 
1 + \frac{x_1+y_1}{2} + &\sum_{n\geq 2}\sum_{i\geq 0}(-1)^i
\ch_{W(B_n)}(H^i(\YBn(\R);\Q))\\
&=(\Sech_B^\sign\circ\Arcsinh^\sign)
(\Exp_B\circ\Arcsinh^\sign)\\
&=\textup{harmonic mean of }
(\Exp\circ\Arcsinh^\sign)|_{p_{2^l}\to x_{2^l}}\\
&\qquad\qquad\textup{ and }
(\Exp\circ\Arcsinh^\sign)|_{p_1\to y_1,\, p_{2^l}\to x_{2^l},l\geq 1}\\
&=\left(\frac{1}{2(1+x_1)}+\frac{1}{2(1+y_1)}\right)^{-1}\circ
(\Exp\circ\Arcsinh^\sign-1).
\end{split}
\ees
\epr
\bpf
The first expression is obtained simply by setting 
$t\to 1$ in Theorem \ref{mainthm}. Now note that
$\Cosh_B^\sign$ is the average of 
\[ \Exp_B^\sign=\exp(\sum_{m\geq 1}\frac{(-1)^{m-1}x_m+(-1)^m y_m}{2m})
\text{ and }\exp(\sum_{m\geq 1}\frac{-x_m+y_m}{2m}), \] 
from which it follows that $\Cosh_B^\sign\Exp_B^{-1}$
is the average of
\[ \exp(-\sum_{\substack{m\geq 2\\m\text{ even}}}\frac{x_m}{m}-
\sum_{\substack{m\geq 1\\m\text{ odd}}}\frac{y_m}{m})\text{ and }
\exp(-\sum_{m\geq 1}\frac{x_m}{m}). \]
So $\Sech_B^\sign \Exp_B$ is the harmonic mean of
$\Exp|_{\substack{p_i\to y_i, i\text{ odd}\\p_i\to x_i, i\text{ even}}}$
and $\Exp|_{p_i\to x_i}$.
From this the second expression follows, and the third expression
is an obvious rewriting of that.
\epf
Although we do not know a general formula for the values of the Euler
character of $\YBn(\R)$, 
it is at least clear from Proposition \ref{eulercharprop}
that it vanishes on elements of $W(B_n)$
which include a cycle of length not a power
of $2$, or a cycle of length $>1$ and cycle-product $-1$.

The reason for including the
third expression in Proposition \ref{eulercharprop} is to point out the
analogy with the result for $\YBn(\C)$, which is the $r=2$ case of
\cite[Theorem 3.4]{hendersonwreath}. Indeed, Proposition \ref{eulercharprop}
can be proved by the same method as the latter result: essentially,
one just uses the fact that Euler characteristic
(in the version with compact supports) is additive
with respect to the nested-set stratification of $\YBn(\R)$. The role of
the series denoted $\CalP(2)$ in \cite{hendersonwreath} is played by
\bes
\begin{split}
\frac{x_1+y_1}{2}+\sum_{n\geq 2}\sum_{i\geq 0}(-1)^i
\ch_{W(B_n)}(&H^i_c(\Pp M_{B_n};\Q))\\
&=1-\frac{1}{2(1+x_1)}-\frac{1}{2(1+y_1)}
\end{split}
\ees
(an easy calculation since $\Pp M_{B_n}$ is homotopy
equivalent to $2^{n-1} n!$ points),
while the role of $\overline{\CalP}(1)$ is played by
$\Exp\circ\Arcsinh^\sign-1$, because of \eqref{typeaeulerchareqn}.
In the complex case one can distinguish
the different cohomology groups by their Hodge weights, even after
taking the alternating sum (the indeterminate $q$ performs this function
in \cite{hendersonwreath}), whereas in the real case this method
gives only the Euler character and not the stronger Theorem \ref{mainthm}.
\section{Moduli space intepretation of $\YBn$}
It was noted in \cite{wonderful} that the variety
$\overline{Y}_{A_n}$ is isomorphic to $\overline{\CalM_{0,n+2}}$,
the moduli space of stable genus $0$ curves
with $n+2$ marked points. (See \cite[Chapter III]{maninbook} for
the definition and basic properties of this variety.)
The map $M_{A_n}\to\overline{Y}_{A_n}$
induces an open embedding $\Pp M_{A_n}\hookrightarrow\overline{Y}_{A_n}$,
whose image corresponds to the open subvariety $\CalM_{0,n+2}$
of $\overline{\CalM_{0,n+2}}$
parametrizing smooth (equivalently, irreducible) such
curves. Explicitly, if we
define the hyperplane arrangement of type $A_n$ in 
$\mathbb{G}_a\!\setminus\!\Aa^{n+1}$, where the additive group 
$\mathbb{G}_a$ acts by simultaneous translation of the coordinates,
then the morphism $M_{A_n}\to \CalM_{0,n+2}$ sends
$(x_1,\cdots,x_{n+1})$ to the point parametrizing $\Pp^1$ with
marked points $x_1,\cdots,x_{n+1},\infty$. 

The variety $\overline{\CalM_{0,n+2}}$ has a well-known stratification,
where the strata are indexed by rooted trees with leaves labelled 
$1,\cdots,n+1$: the open stratum $\CalM_{0,n+2}$ corresponds to the
tree whose only vertices are the root and the leaves, and more generally
the other vertices of the tree specify the pattern of intersections
of the components of the stable curve (the root corresponds to the
last marked point). This stratification corresponds in a natural
way to the nested-set stratification of $\overline{Y}_{A_n}$ 
defined in \cite{wonderful}.

One striking consequence of the isomorphism $\overline{Y}_{A_n}\cong
\overline{\CalM_{0,n+2}}$ is that the action of the
linear automorphism group $W(A_n)=S_{n+1}$ can be extended to
an action of $S_{n+2}$ (permuting the marked points). The action of
the `non-linear' symmetries on cohomology
is difficult to describe from the point of view of Theorem
\ref{rainscoxeterthm}. Nevertheless, in \cite{rainstypea} the 
formula for the character of $S_{n+1}$ on the rational cohomology
groups of $\overline{\CalM_{0,n+2}}(\R)$ resulting from
Theorem \ref{rainstypeathm} was extended to all of $S_{n+2}$.

Now we want to give a similar moduli space interpretation of
$\YBn$. Realize the hyperplane complement $M_{B_n}$ in the usual way:
\[ M_{B_n}=\{(x_1,\cdots,x_n)\in\Aa^n\,|\,x_i\neq 0\text{ for all }i,\
x_i\neq \pm x_j\text{ for }i\neq j\}. \]
We can then define a closed embedding $\pi:M_{B_n}\to M_{A_{2n}}$ by
\beq 
\pi(x_1,\cdots,x_n)=(x_1,\cdots,x_n,-x_1,\cdots,-x_n,0),
\eeq
where the right-hand side represents a $\mathbb{G}_a$-orbit.
Clearly this embedding descends to an embedding 
$\pi:\Pp M_{B_n}\hookrightarrow \Pp M_{A_{2n}}\cong \CalM_{0,2n+2}$,
whose image $\pi(\Pp M_{B_n})$ is the closed subvariety parametrizing
smooth genus $0$ curves
with $2n+2$ marked points and an involution which fixes
the last two marked points and interchanges the $j$th and $(n+j)$th
for $1\leq j\leq n$ (such an involution is clearly unique if it exists).
\bth \label{modulithm}
There is a closed embedding $\pi:\YBn\hookrightarrow\overline{\CalM_{0,2n+2}}$
whose restriction to $\Pp M_{B_n}$ is the embedding
$\Pp M_{B_n}\hookrightarrow \CalM_{0,2n+2}$ described above. The image
$\pi(\YBn)$ is the subvariety parametrizing stable genus $0$ curves
with $2n+2$ marked points and an involution which fixes
the last two marked points, interchanges the $j$th and $(n+j)$th
for $1\leq j\leq n$, and does not fix any component pointwise
(such an involution is unique if it exists). The stratification
of $\YBn$ obtained by intersecting this image with the strata of
$\overline{\CalM_{0,2n+2}}$ is the nested-set stratification of 
\cite{wonderful}.
\eth
\bpf
The first statement follows 
from the fact that $\pi$ induces a surjective map of
building sets $\CalG_{A_{2n}}\to\CalG_{B_n}$; see \cite[Proposition 2.4 and
the last statement of Corollary 2.6]{rainshomology}.
 
The image $\pi(\YBn)$ equals the closure in $\overline{\CalM_{0,2n+2}}$
of the subvariety $\pi(\Pp M_{B_n})\subseteq \CalM_{0,2n+2}$
described above; we must show that this equals the
subvariety 
$X\subseteq\overline{\CalM_{0,2n+2}}$
parametrizing stable curves with an involution as stated.
That $X$ is contained in the closure of $\pi(\Pp M_{B_n})$ follows
from the fact that every stable genus $0$ curve with $2n+2$
marked points and an involution of the stated form can be obtained
by progressively degenerating a smooth one (maintaining
the existence of an involution of the stated form throughout the
degeneration). For example:
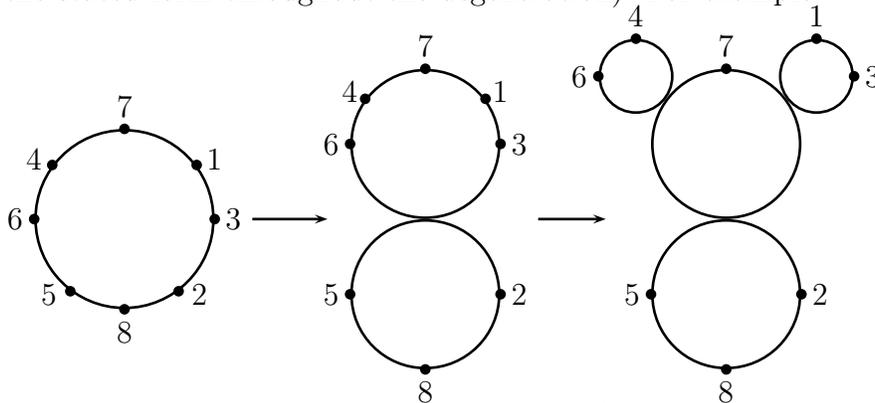
\begin{center}
\begin{pspicture}(5.55,-2)(16.95,2.5)
    	\psset{linewidth=1pt}%
	\pscircle(7,0){1.2}
	\psdot(7,1.2)
	\rput(7,1.5){$7$}
	\psdot(7,-1.2)
	\rput(7,-1.5){$8$}
	\psdot(8.2,0)
	\rput(8.45,0){$3$}
	\psdot(5.8,0)
	\rput(5.55,0){$6$}
	\psdot(7.96,0.72)
	\rput(8.2,0.8){$1$}
	\psdot(6.04,0.72)
	\rput(5.8,0.8){$4$}
	\psdot(7.72,-0.96)
	\rput(8,-1){$2$}
	\psdot(6.28,-0.96)
	\rput(6,-1){$5$}
	\psline{->}(8.7,0)(9.7,0)
	\pscircle(11,1){1}
	\psdot(11,2)
	\rput(11,2.3){$7$}
	\psdot(11.8,1.6)
	\rput(12,1.7){$1$}
	\psdot(10.2,1.6)
	\rput(10,1.7){$4$}
	\psdot(12,1)
	\rput(12.25,1){$3$}
	\psdot(10,1)
	\rput(9.75,1){$6$}
	\pscircle(11,-1){1}
	\psdot(11,-2)
	\rput(11,-2.3){$8$}
	\psdot(12,-1)
	\rput(12.25,-1){$2$}
	\psdot(10,-1)
	\rput(9.75,-1){$5$}
	\psline{->}(12.5,0)(13.4,0)
	\pscircle(15,1){1}
	\psdot(15,2)
	\rput(15,2.3){$7$}
	\pscircle(15,-1){1}
	\psdot(15,-2)
	\rput(15,-2.3){$8$}
	\psdot(16,-1)
	\rput(16.25,-1){$2$}
	\psdot(14,-1)
	\rput(13.75,-1){$5$}
	\pscircle(16.2,1.9){0.5}
	\psdot(16.2,2.4)
	\rput(16.2,2.7){$1$}
	\psdot(16.7,1.9)
	\rput(16.95,1.9){$3$}
	\pscircle(13.8,1.9){0.5}
	\psdot(13.8,2.4)
	\rput(13.8,2.7){$4$}
	\psdot(13.3,1.9)
	\rput(13.05,1.9){$6$}
\end{pspicture}
\end{center}
Here circles represent irreducible components (whose intersections are
transverse, despite appearances), and arrows indicate a passage to the
limit in some family of curves; the claim is just that the vertical axis
of symmetry can be maintained throughout.

Now take any point $x_C$
in the closure of $\pi(\Pp M_{B_n})$, corresponding to the isomorphism
class of a stable genus $0$ curve $C$ with $2n+2$ marked points; since
the existence of an involution fixing
the last two marked points and interchanging the $j$th and $(n+j)$th
for $1\leq j\leq n$ is a closed property, $C$ has such an involution, say
$\tau$. We must show that $\tau$ does not fix any component
of $C$ pointwise. Suppose that it did; it is easy to see that there
must be some pointwise-fixed component $C_1$ such that the point
in $C_1$ `closest to' the $(2n+1)$th marked point is different from
the point in $C_1$ `closest to' the $(2n+2)$th marked point. By
stability, $C_1$ must have at least one other special point, which is
the point in $C_1$ `closest to' the $j$th marked point and hence also
the $(n+j)$th marked point for some $1\leq j\leq n$. For example, suppose 
$n=3$ and $C$ is:
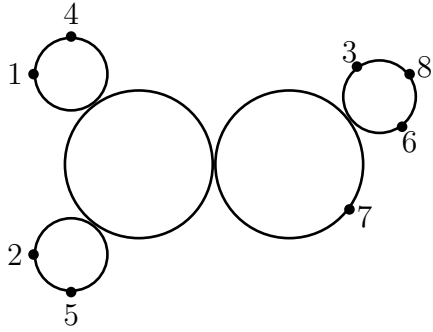
\begin{center}
\begin{pspicture}(1.35,-1.8)(6.7,1.8)
   	\psset{linewidth=1pt}%
	\pscircle(3,0){1}
	\pscircle(2.1,1.2){0.5}
	\psdot(1.6,1.2)
	\rput(1.35,1.2){$1$}
	\psdot(2.1,1.7)
	\rput(2.1,2){$4$}
	\pscircle(2.1,-1.2){0.5}
	\psdot(1.6,-1.2)
	\rput(1.35,-1.2){$2$}
	\psdot(2.1,-1.7)
	\rput(2.1,-2){$5$}
	\pscircle(5,0){1}
	\psdot(5.8,-0.6)
	\rput(6,-0.7){$7$}
	\pscircle(6.2,0.9){0.5}
	\psdot(5.9,1.3)
	\rput(5.8,1.5){$3$}
	\psdot(6.6,1.2)
	\rput(6.8,1.3){$8$}
	\psdot(6.5,0.5)
	\rput(6.6,0.3){$6$}
\end{pspicture}
\end{center}
Here all components are fixed by $\tau$, but only the two represented
by larger circles are pointwise-fixed. Of these, only the right-hand one
satisfies the condition of $C_1$. We can take
either $j=1$ or $j=2$, but not $j=3$.

Now there is a projective morphism
$p_j:\overline{\CalM_{0,2n+2}}\to\overline{\CalM_{0,4}}\cong\Pp^1$ defined by
`stably forgetting' all marked points except the $j$th, $(n+j)$th,
$(2n+1)$th and $(2n+2)$th, and what we know about $C_1$ implies that
$x_C$ belongs to the
preimage $p_j^{-1}(\{b\})$, where $b$ is one of the three points of
$\overline{\CalM_{0,4}}\setminus \CalM_{0,4}$. On the other hand,
$\pi(\Pp M_{B_n})$ is contained in the
preimage $p_j^{-1}(\{a\})$ where $a\in \CalM_{0,4}$
parametrizes the curve $\Pp^1$ with marked points
$1,-1,0,\infty$. This contradicts the assumption that $x_C$ was
in the closure of $\pi(\Pp M_{B_n})$.

Finally, it is easy to see that a
stratum of $\overline{\CalM_{0,2n+2}}$ intersects the subvariety
$\pi(\YBn)$ if and only if the corresponding tree admits an involution
which fixes the root and the $(2n+1)$th leaf, interchanges the
$j$th and $(n+j)$th leaves for $1\leq j\leq n$, and does not
fix any edges except those forming the path from the root to the
$(2n+1)$th leaf (such an involution is unique if it exists).
This set of trees can be easily identified with the set 
denoted $\CalT(2,n)$ in \cite[Section 6]{hendersonwreath}, and
the final statement follows from what was said there.
\epf
One corollary is that the sequence of varieties $(\YBn)$ form a right module
for the operad $(\overline{\CalM_{0,n}})$; for $1\leq j\leq n$, the operation
\[ \circ_j:\YBn\times\overline{\CalM_{0,m+1}}\to\overline{Y}_{B_{n+m-1}} \]
is defined by gluing two copies of the second stable curve to the first,
at the $j$th and $(n+j)$th marked points, and renumbering marked points
appropriately. A similar construction is made in \cite{armstrongetal}.

Observe that if $C$ is a stable genus $0$ curve with $2n+2$ marked
points and an involution $\tau$ as in Theorem \ref{modulithm},
then the fixed points of $\tau$ are exactly the last two marked
points and the finite set of nodes which lie on the `path' between them.
The last marked point is thus determined by the rest of the data,
and $C$ is still stable without it, so $\YBn$ can also be identified
with a subvariety of $\overline{\CalM_{0,2n+1}}$. Indeed, the above proof
would have worked just as well for the embedding 
$M_{B_n}\hookrightarrow M_{A_{2n-1}}$ obtained by removing the last
coordinate of $\pi:M_{B_n}\hookrightarrow M_{A_{2n}}$.

The advantage of the description in Theorem \ref{modulithm} is that it
makes clear some extra symmetry. The subgroup of $S_{2n+2}$ which
preserves $\pi(\YBn)$ is the product $W(B_n)\times\langle(2n+1,2n+2)\rangle$,
where $W(B_n)$ is embedded in $S_{2n}$ as the centralizer of 
$j\leftrightarrow n+j$. The resulting action of $W(B_n)$ on $\YBn$ is the one
we have been considering throughout the paper; the action of
$(2n+1,2n+2)$ gives an involution $\sigma$ of $\YBn$ which commutes
with the $W(B_n)$-action.
Note that the restriction of $\sigma$ to $\Pp M_{B_n}$
is the map $(x_1,\cdots,x_n)\mapsto (x_1^{-1},\cdots,x_n^{-1})$.
\bpr
The fixed subvariety $(\YBn)^\sigma$ is isomorphic to the disjoint
union of $2^{n-1}$ copies of $\overline{\CalM_{0,n+1}}$.
\epr
\bpf
Suppose that $x_C\in \pi(\YBn)$ corresponds to the isomorphism
class of a stable curve $C$ with involution $\tau$ as above. Then
$\sigma(x_C)=x_C$ if and only if
$C$ has another involution $\tau'$ which fixes the first
$2n$ marked points and interchanges the last two. This implies that
the last two marked points lie on the same component $C_1$, since
otherwise there would be another marked point which was `closer to'
one than to the other. It is also clear that $C_1$ has exactly two
other special points $a$ and $\tau(a)$, which are the fixed points of 
$\tau'$ in $C_1$. If $C'$ denotes the union of the components $\neq C_1$
to which $a$ is the closest point of $C_1$, then $C'$ is a stable
curve with $n+1$ marked points (the $n$ of the original $2n+2$
which it contains, together with $a$). Since the division of the
first $2n$ marked points into two parts interchanged by $\tau$
can be made in $2^{n-1}$ ways, the result follows.
\epf
\bcr \label{sigmalefschetzcor}
The Lefschetz number $\sum_i (-1)^i \tr(\sigma,H^i(\YBn(\R);\Q))$ equals
$0$ if $n$ is odd, and 
\[ \frac{(-1)^\frac{n-2}{2}n!(n-2)!}{(\frac{n}{2})!(\frac{n-2}{2})!}
\text{ if $n$ is even.} \]
\ecr
\bpf
By a well-known principle (see \cite{wall}, for instance),
the Lefschetz number is $\chi(\YBn(\R)^\sigma)$, which equals
$2^{n-1}\chi(\overline{Y}_{A_{n-1}}(\R))$ by
the Proposition.
Thus \eqref{typeaeulercharacteristiceqn} gives the result.
\epf
One could perhaps determine
$\sum_i (-1)^i \tr(w\sigma,H^i(\YBn(\R);\Q))$ for general $w\in W(B_n)$
by similar arguments.

We can describe the action of $\sigma$ on $H^1(\YBn(\R);\Q)$ using an
explicit basis. Theorem \ref{rainscoxeterthm} showed that
$H^1(\YBn(\R);\Q)$ is canonically the direct sum
of one-dimensional subspaces corresponding to the irreducible rank-$2$
parabolic subgroups of $W(B_n)$, or equivalently to those factors
in the right-hand side of \eqref{modeldefeqn} which are isomorphic to
$\Pp^1$. More precisely, one can 
construct a basis using the projections to these 
factors, as follows (recall the notation for the building set $\CalG_{B_n}$
introduced in the proof of Proposition \ref{typedprop}). 
For any $i\neq j$ in $\{1,\cdots,n\}$, define
a morphism $\rho_{ij}:\YBn\to\Pp^1$ as the composition of the projection
$\YBn\to\Pp(V/(\R^{\{i,j\}})^\perp)$ and the identification
$\Pp(V/(\R^{\{i,j\}})^\perp)\cong\Pp^1$ which takes $x_i$ and 
$x_j$ as the
homogeneous coordinates, in that order. Thus the restriction of
$\rho_{ij}$ to $\Pp M_{B_n}$ is the map 
$(x_1,\cdots,x_n)\mapsto x_i x_j^{-1}$.
Denote by $\nu_{ij}$ the pull-back 
$\rho_{ij}^*(\alpha)\in H^1(\YBn(\R);\Q)$,
where $\alpha\in H^1(\Pp^1(\R);\Q)$ is the standard class.
For any distinct $i,j,k$ in $\{1,\cdots,n\}$ and 
$\varphi:\{i,j,k\}\to\{\pm 1\}$, define 
$\rho_{ijk,\varphi}:\YBn\to\Pp^1$ as the composition of the projection
$\YBn\to\Pp(V/(\R_\varphi^{\{i,j,k\}})^\perp)$ and the identification
$\Pp(V/(\R_\varphi^{\{i,j,k\}})^\perp)\cong\Pp^1$ which 
takes $\varphi(i)x_i-\varphi(k)x_k$ and 
$\varphi(j)x_j-\varphi(k)x_k$ as the homogeneous coordinates, in 
that order. Thus the restriction of
$\rho_{ijk,\varphi}$ to $\Pp M_{B_n}$ is the map 
$(x_1,\cdots,x_n)\mapsto (\varphi(i)x_i-\varphi(k)x_k)
(\varphi(j)x_j-\varphi(k)x_k)^{-1}$. Denote by
$\omega_{ijk,\varphi}$ the pull-back
$\rho_{ijk,\varphi}^*(\alpha)\in H^1(\YBn(\R);\Q)$. 
\bpr 
\ben
\item $\nu_{ji}=-\nu_{ij}$, $\omega_{ijk,-\varphi}=\omega_{ijk,\varphi}$, and
\[ \omega_{i_{w(1)}i_{w(2)}i_{w(3)},\varphi}=\varepsilon(w)
\omega_{i_1 i_2 i_3,\varphi}
\text{ for $w\in S_3$.} \]
\item For any element $(\varepsilon_1,\cdots,\varepsilon_n)w$ of $W(B_n)$,
where $(\varepsilon_1,\cdots,\varepsilon_n)\in\Gamma_n$ and $w\in S_n$,
\bes 
\begin{split}
(\varepsilon_1,\cdots,\varepsilon_n)w.\nu_{ij}
&=\varepsilon_{w(i)}\varepsilon_{w(j)}\nu_{w(i)w(j)},\\
(\varepsilon_1,\cdots,\varepsilon_n)w.\omega_{ijk,\varphi}
&=\omega_{w(i)w(j)w(k),(w(a)\mapsto \varepsilon_{w(a)}\varphi(a))}.
\end{split}
\ees
\item The disjoint union of $\{\nu_{ij}\,|\,1\leq i<j\leq n\}$ and
$\{\omega_{ijk,\varphi}\,|\,1\leq i<j<k\leq n\}$, where to avoid
redundancy we take only one from each pair $\pm\varphi$, is a basis
of $H^1(\YBn(\R);\Q)$.
\item The span of $\{\omega_{ijk,\varphi}\}$
is the image of the natural injective map 
$H^1(\YDn(\R);\Q)\to H^1(\YBn(\R);\Q)$
mentioned in Section 1.
\item The span of $\{\nu_{ij}\}$
is the $(-1)$-eigenspace of $\sigma$ on $H^1(\YBn(\R);\Q)$.
\item $\sigma(\omega_{ijk,\varphi})=\omega_{ijk,\varphi}-
\varphi(i)\varphi(j)\nu_{ij}-\varphi(j)\varphi(k)\nu_{jk}-
\varphi(k)\varphi(i)\nu_{ki}$.
\een
\epr
\bpf
Parts (1) and (2) follow easily from the definitions, and parts 
(3) and (4) from
the construction of the isomorphism in \cite[Theorem 3.7]{rainshomology}.
For parts (5) and (6), we use the embedding 
$\pi:\YBn\hookrightarrow\overline{\CalM_{0,2n+2}}$. Recall the elements
$\omega_{ijkl}\in H^1(\overline{\CalM_{0,2n+2}}(\R);\Q)$ from
\cite{etingofetal}, defined as $\rho_{ijkl}^*(\alpha)$
for certain morphisms $\rho_{ijkl}:\overline{\CalM_{0,2n+2}}\to\Pp^1$,
where $i,j,k,l$ are distinct elements of $\{1,\cdots,2n+2\}$.
It is clear from the definitions that the morphism $\rho_{ij}$ used above
is the composition $\rho_{i,j,2n+1,2n+2}\circ \pi$, and the morphism
$\rho_{ijk,\varphi}$ is the composition $\rho_{i^\varphi,j^\varphi,k^\varphi,
2n+2}\circ \pi$, where $i^\varphi$ denotes $i$ if $\varphi(i)=1$ and
$i+n$ if $\varphi(i)=-1$. Hence
\bes
\nu_{ij}=\pi^*(\omega_{i,j,2n+1,2n+2}),\
\omega_{ijk,\varphi}=\pi^*(\omega_{i^\varphi,j^\varphi,k^\varphi,2n+2}).
\ees
Using the definition of $\sigma$ and the relation \cite[(2.1)]{etingofetal},
we calculate
\bes
\begin{split}
\sigma(\nu_{ij})&=\pi^*(\omega_{i,j,2n+2,2n+1})=\pi^*(-\omega_{i,j,2n+1,2n+2})
=-\nu_{ij},\\
\sigma(\omega_{ijk,\varphi})&=\pi^*(\omega_{i^\varphi,
j^\varphi,k^\varphi,2n+1})\\
&=\pi^*(\omega_{i^\varphi,
j^\varphi,k^\varphi,2n+2}\\
&\qquad
-\omega_{i^\varphi,j^\varphi,2n+1,2n+2}-\omega_{j^\varphi,k^\varphi,2n+1,2n+2}
-\omega_{k^\varphi,i^\varphi,2n+1,2n+2})\\
&=\omega_{ijk,\varphi}-
\varphi(i)\varphi(j)\nu_{ij}-\varphi(j)\varphi(k)\nu_{jk}-
\varphi(k)\varphi(i)\nu_{ki}.
\end{split}
\ees
Thus part (6) is proved, and the span of $\{\nu_{ij}\}$ is at least contained
in the $(-1)$-eigenspace of $\sigma$; however, (6) shows that $\sigma$
acts trivially on the quotient by the span of $\{\nu_{ij}\}$, which proves (5).
\epf
The fixed point subspace $H^1(\YBn(\R);\Q)^\sigma$ is spanned by the elements
\[ \widetilde{\omega}_{ijk,\varphi}=\omega_{ijk,\varphi}-
\frac{1}{2}\varphi(i)\varphi(j)\nu_{ij}-\frac{1}{2}\varphi(j)\varphi(k)\nu_{jk}
-\frac{1}{2}\varphi(k)\varphi(i)\nu_{ki}, \]
and is isomorphic to $H^1(\YDn(\R);\Q)$ as a representation of $W(B_n)$
via the map $\widetilde{\omega}_{ijk,\varphi}\mapsto\omega_{ijk,\varphi}$.
However, this isomorphism does not generalize to higher-degree cohomologies
(indeed, Corollary \ref{sigmalefschetzcor} shows that
$\dim H^2(\overline{Y}_{B_4}(\R);\Q)^\sigma=3$, whereas 
$\dim H^2(\overline{Y}_{D_4}(\R);\Q)=15$).

We conjecture that the cohomology ring 
$H^\bullet(\overline{Y}_{B_n}(\R);\Q)$ is generated in degree $1$. Since the
analogous statement for $\overline{M_{0,n}}(\R)$ was proved in
\cite{etingofetal}, it is equivalent to conjecture that the homomorphism
$\pi^*:H^\bullet(\overline{M_{0,2n+2}}(\R);\Q)\to 
H^\bullet(\overline{Y}_{B_n}(\R);\Q)$ is surjective. Calculations
for small $n$ suggest that the only relations satisfied by the
above basis elements of $H^1(\overline{Y}_{B_n}(\R);\Q)$ are those
forced by the relations in \cite[Definition 2.1]{etingofetal};
more precisely, that the ideal 
$\ker(\pi^*)$ of $H^\bullet(\overline{M_{0,2n+2}}(\R);\Q)$ 
is also generated in degree $1$.

\end{document}